\documentclass{amsart}

\usepackage{amssymb,mathrsfs,amsmath,dsfont,array}
\usepackage{multirow,arydshln}
\usepackage{tikz}
\usetikzlibrary{arrows,positioning,decorations.pathmorphing,shapes}
\usepackage{capt-of}
\usepackage{mathbbol}     
\usepackage{color}
\usepackage{verbatim} 
\usepackage[all,cmtip]{xy}
\usepackage{amsbsy}
\usepackage{amsmath}
\usepackage{nccmath}

\usepackage{hyperref}

\newtheorem{Thm}{Theorem}[section]

\newtheorem{Pro}[Thm]{Proposition}

\theoremstyle{definition}
\newtheorem{deft}[Thm]{Definition}
\newtheorem{exa}[Thm]{Example}

\theoremstyle{remark}
\newtheorem{rem}[Thm]{Remark}

\numberwithin{equation}{section}

\def\a1s{a_1,\cdots, a_s}
\def\a{\alpha}

\def\aa{\mathcal A}
\def\aA{\mathscr{A}}

\def\andd{\quad\hbox{and}\quad}

\def\b{\beta}

\def\bs{\boldsymbol}

\def\bl4{B_{\ell\geq4}}

\def\bbbc{{\mathbb C}}

\def\d{\delta}
\def\D{\Delta}

\def\fg{\mathfrak{g}}

\def\scg{\mathscr{G}}

\def\hB{\mathscr{B}}

\def\fh{\mathfrak{h}}

\def\fk{\mathfrak{k}}

\def\lam{\lambda}
\def\Lam{\Lambda}
\def\LL{\mathcal{L}}

\def\fl{\mathfrak{L}}

\def\ep{\epsilon}
\def\fm{(\cdot,\cdot)}

\def\bbbq{\mathbb{Q}}

\def\bbbr{{\mathbb R}}

\def\supp{\hbox{\rm supp}}

\def\1k{\frac{1}{k}}
\def\op{\oplus}
\def\ot{\otimes}

\def\sub{\subseteq}
\def\sg{\sigma}

\def\pf{\noindent{\bf Proof. }}

\def\sspan{\hbox{\rm span}}

\def\fp{\mathfrak{p}}

\def\ft{\mathfrak{t}}

\def\w{{\mathcal W}}

\def\bbbz{{\mathbb Z}}

\def\1il{1\leq i\leq\ell}

\begin{document}

 \title[Quasi-Integrable Modules]{Quasi-Integrable Modules,\\ a Class of Non-highest Weight Modules  over Twisted Affine Lie Superalgebras}


\author{Malihe Yousofzadeh}
\address{ Department of Pure Mathematics, Faculty of Mathematics and Statistics, University of Isfahan, Isfahan, Iran, P.O.Box 81746-73441 and School of Mathematics, Institute for Research in
Fundamental Sciences (IPM), P.O. Box: 19395-5746, Tehran, Iran.
 }
\curraddr{}
\email{ma.yousofzadeh@sci.ui.ac.ir \& ma.yousofzadeh@ipm.ir.}
\thanks{This research  was in part supported by a grant from IPM (No. 1400170215) and  is partially carried out in
IPM-Isfahan Branch.}


\subjclass[2020]{17B10, 17B67\\ Keywords: Twisted Affine Lie Superalgebras, Finite Weight Module. }

\date{}

\dedicatory{}

\commby{}

\begin{abstract}
In this paper, we characterize quasi-integrable modules, of nonzero level, over twisted affine Lie superalgebras. We show that quasi-integrable modules are not necessarily highest weight modules. We prove that each quasi-integrable module  is parabolically induced from  a cuspidal module, over a finite dimensional  Lie superalgebra having a Cartan subalgebra whose corresponding root system just contain real roots; in particular, the classification of quasi-integrable modules is reduced to the known classification of cuspidal modules over such Lie superalgebras.
\end{abstract}
\maketitle

\section{Introduction}
Affine Lie superalgebras,  which  are the super version of affine Lie algebras, were introduced and classified by Van de Leur in 1986. Affine Lie superalgebras are divided into (1) nontwisted types $X^{(1)},$ where $X$ is the type of a finite dimensional basic classical simple Lie superalgebra, and (2) twisted types $A(2k-1,2\ell-1)^{(2)}$  ({\tiny$(k,\ell)\neq (1,1)$}), $A(2k,2\ell)^{(4)},$  $A(2k,2\ell-1)^{(2)}$ and $D(k,\ell)^{(2)}$ together with twisted affine Lie algebras.  Roughly speaking, a nontwisted affine Lie superalgebra is a certain  extension of a current Lie superalgebra $\fk\ot_\bbbc \bbbc[t^{\pm1}]$ in which $\fk$ is a basic classical simple Lie superalgebra but to construct a twisted affine Lie superalgebra, a finite order automorphism  gets involved as well.

Representation theory of affine Lie (super)algebras is one of the important topics which has drawn  considerable attention of mathematicians as well as physicists. In this regard, after the  study of finite dimensional modules, the first step is the study of finite weight modules.  In 1986 and 1988, V.~Chari and A.~Perssley (\cite{C} \& \cite{CP2})  studied  integrable finite weight modules over an affine Lie algebra and showed that the only irreducible integrable finite weight modules of nonzero level are    irreducible integrable highest weight modules and their duals.  Since then, there have been several attempts to study representation theory of affine Lie superalgebras; see  \cite{DG},  \cite{E1}, \cite{E2},  \cite{F5}, \cite{FS} and the references therein.

For most  affine Lie superalgebras with nonzero odd part, the even part contains  two affine Lie subalgebras in twisted case and  two or three affine Lie subalgebras in nontwisted case.  In \cite{E1}, the author shows that if $M$ is an   nonzero level irreducible integrable finite weight module over a  nontwisted affine Lie superalgebra, then $M$ is a trivial module if the even part  has at least two affine components and it is a highest weight module if the even part has only one affine component; see also \cite{KW}.

The structure of an irreducible  finite weight module $M$  over a  twisted affine Lie superalgebra $\fl$ with root system $R$, strongly depends on the nature of the action of root vectors corresponding to nonzero real roots, i.e., roots which are not self-orthogonal. More precisely, each nonzero root vector corresponding to a nonzero real root $\a$, acts on $M$ either injectively or locally nilpotently. We denote by $R^{in}$ (resp. $R^{ln}$) the set of all nonzero real roots $\a$ whose corresponding  nonzero root vectors  act on $M$  injectively (resp.  locally nilpotently).  In \cite{you8}, we showed that if $\{\d\}$ generates the group generated by the set $R_{im}$ consisting of all roots which are orthogonal to all other roots, for each nonzero real root $\a,$ one of the following occurs:

\noindent$\bullet$~  $\a$ is full-locally nilpotent, i.e., $R\cap (\a+\bbbz\d)\sub R^{ln},$

\noindent$\bullet$~ $\a$ is full-injective, i.e.,  $R\cap (\a+\bbbz\d)\sub R^{in},$

\noindent$\bullet$~$\pm\a$ are up-nilpotent (resp. down-nilpotent) hybrid, i.e., there is  $m\in\bbbz^{>0}$ with
$$\hbox{\small $R\cap (\pm\a+\bbbz^{\geq m}\d)\sub R^{ln}$ (resp.  $R^{in}$) and $R\cap (\pm\a+\bbbz^{\leq -m}\d)\sub R^{in}$ (resp. $R^{ln}$)}.$$
So,  we can divide our study into two cases either all nonzero real roots are hybrid or some of nonzero real roots are not hybrid; we call $M$ respectively hybrid or tight.

In \cite{you8}, we reduced  the classification of  hybrid  irreducible finite weight modules over  a twisted affine Lie superalgebra to the classification of  cuspidal modules of finite dimensional  Lie superalgebras; see \cite{DMP}.
 In \cite{you9}, we began the study of tight modules over a twisted  affine Lie superalgebra $\fl$ whose even part has two affine components.
  The structure of such modules  depends on  whether non-hybrid roots occur for the roots of both affine components  or not.
  We called the root systems of affine components of $\fl_0,$ $R(1)$ and $R(2)$ and showed that if for a tight module with bounded weight multiplicity, $R(i)\cap R^{ln}$ ($i=1,2$) is a proper nonempty subset of the set of real roots of $R(i),$ then $M$ is parabolically induced.
In the present paper, we study the case that for a choice $i,j$ with $\{i,j\}=\{1,2\},$ all nonzero real roots of $R(i)$ belong to $R^{ln}$ while all nonzero real roots of $R(j)$ are hybrid; we  call such modules quasi-integrable modules.
 We also give an example of a module which is quasi-integrable while it is not a highest weight module.

This work is some part of a long-term project on representation theory of twisted affine Lie superalgebras. The main aim of this part is the classification of quasi-integrable modules over twisted affine Lie superalgebras. This in turn is divided into two subparts. In the  first subpart,  we give the structure of quasi-integrable modules and in the second part, we give the classification. In this paper, we focus on the first subpart.
We prove that each quasi-integrable irreducible finite weight  module is parabolically induced from  a cuspidal module, over a direct sum of a reductive Lie algebra and finitely many basic classical simple Lie superalgebras whose root systems have just real roots. This is the most important step as the classification of cuspidal modules over such Lie superalgebras are known by \cite{M}, \cite{FGG} and \cite{G}.

\section{Some notations and definitions}\label{preliminaries}
In  this section, to have a self-contained paper,  we gather some information we need throughout the paper.
 All vector spaces and tensor products are considered over the field of complex numbers $\bbbc$. The degree of a homogeneous element $x$ of a superspace is denoted by $|x|.$ We mention that whenever we use $|x|$ for an element $x$ of a superspace, we mean that $x$ is homogeneous.

Suppose that  $\fk=\fk_0\op\fk_1$  is  a Lie superalgebra, $\fh$ is a finite dimensional subalgebra of the even part $\fk_0$ of $\fk$ and $\fm$ is  a   nondegenerate supersymmetric invariant even bilinear form on $\fk.$
 A superspace  ${V}={V}_0\op {V}_1$ is called  a {\it weight  $\fk$-module} (with respect to $\fh$) if
 $[x,y]v=x(yv)-(-1)^{|x||y|}y(xv)$ ($i,j\in\{0,1\}$, $v\in {V},$  $x,y\in\fk$), $\fk^i {V}^j\sub {V}^{i+j}$
 and  ${V}=\op_{\lam\in \fh^\ast}{V}^\lam$ with
${V}^\lam:=\{v\in {V}\mid hv=\lam(h)v\;\;(h\in\fh)\}$ for each $\lam\in \fh^\ast.$ In this case, an element $\lam$ of the {\it support} $\supp({V}):=\{\lam\in \fh^*\mid {V}^\lam\neq \{0\}\}$  of $V$, is called a {\it weight} of $V$ and the corresponding ${V}^\lam$ is called a {\it weight space}. Elements of a weight space are called {\it weight vectors}. If all weight spaces are  finite dimensional, the module ${V}$ is called a {\it finite weight module} (with respect to $\fh$).

If $\fk$ has a weight space decomposition with respect to $\fh$ via the adjoint representation, we say $\fk$  has   a root space decomposition with respect to $\fh$; the set of weights of $\fk$ is called the {\it root system} and weights, weight vectors and weight spaces are called respectively {\it roots}, {\it root vectors} and  {\it root spaces}.

If  $\fk$ has  a root space decomposition with respect to $\fh$ with root system ${\D}$ such that   $\fk^0=\fh,$  the  form on $\fk$ is restricted to  a nondegenerate symmetric bilinear form on $\fh$ which in turn induces naturally   a symmetric nondegenerate bilinear form, denoted again by $\fm$, on the dual space $\fh^*$ of  $\fh$.
So, we can define
 \begin{equation}
  \label{decom}
  \begin{array}{lll}
\hbox{\footnotesize $\D_{im} =\{\a\in \D\mid(\a,\D)=\{0\}\},$}&
\hbox{\footnotesize $\D ^\times_{re} = \{\alpha \in \D \mid (\alpha,\alpha) \neq 0\},$}&
\hbox{\footnotesize $\D _{re} = \D ^\times_{re} \cup \{0\},$}\\
\hbox{\footnotesize $\D ^\times = \D \setminus \D_{im} ,$}&
\hbox{\footnotesize $\D ^\times_{ns} =\{\alpha \in \D^\times | (\alpha,\alpha) = 0\},$}&
\hbox{\footnotesize $\D_{ ns} = \D^ \times_{ns} \cup \{0\}.$}
\end{array}
\end{equation}
Elements of $\D_{im}$ (resp. $\D_{re}$ and $\D_{ns}$) are called imaginary roots (resp. real roots and nonsingular roots).
 We also  have
\begin{equation}
\label{talpha}
\parbox{4.5in}{for each $\a\in \fh^*,$ there is a unique $t_\a\in\fh$ with $\a(h)=(t_\a,h)$ for all $h\in\fh.$}
\end{equation}
In this case, if  ${V}$ is a $\fk$-module and $S\sub {\D},$
\begin{equation}\label{in-nil}
\parbox{4.5in}{we denote by  ${S}^{ln}({V})$ (resp. ${S}^{in}({V})$),   the set of all roots $\a\in S$ with $(\a,\a)\neq 0$ for which $0\neq x\in\fk^\a$ acts on ${V}$ locally nilpotently (resp. injectively).}
\end{equation}
If moreover, $V$ has a  weight space decomposition $V=\op_{\lam\in \fh^*}V^\lam$ with respect to $\fh$ and  $W$ is an $\fh$-submodule of $V,$  then, we have    $W=\op_{\lam\in\fh^*}(W\cap V^\lam);$ we next set\footnote{We use $\#$ to indicate the cardinal number.}
\small{\begin{equation}\label{BC}
\begin{array}{l}
\frak{B}_W:=\{\a\in \hbox{span}_\bbbz {\D}\mid  \#\{k\in\bbbz^{>0}\mid \lam+k\a\in \supp(W)\}<\infty \;(\forall\lam\in\supp(W))\},\\
\frak{C}_W:=\{\a\in \hbox{span}_\bbbz {\D}\mid   \a+\supp(W)\sub\supp(W)\}.
\end{array}
\end{equation}}
We say the $\fk$-module $V$  to have {\it shadow} if
 \begin{itemize}
 \item[\rm\bf (s1)]
${\D}_{re}^\times=\{\a\in {\D}\mid (\a,\a)\neq 0\}={\D}^{in}(V)\cup {\D}^{ln}(V),$

\item[\rm\bf (s2)] ${\D}^{ln}(V)=\frak{B}_V\cap {\D}_{re}^\times$ and ${\D}^{in}(V)=\frak{C}_V\cap {\D}_{re}^\times.$
 \end{itemize}
Suppose that  $V$ is a $\fk$-module having shadow and $\scg$ is a subalgebra of $\fk$ containing $\fh$. Denote the root system of $\scg$ with respect to $\fh$ by $T$. If $W$ is a $\scg$-submodule of the $\scg$-module $V$, one can easily see that $W$ has also shadow and that  $R^{*}(V)\cap T=T^*(W)$ for $*=in,ln$.  In other words, that root vectors corresponding to real roots act locally nilpotently or injectively, depends only on $V$, so if $V$ is fixed, we may simply  denote  $R^*(V)$ and $T^*(W)$ by $R^*$ and $T^*$
respectively, for $*=in,ln.$

\smallskip

For a subset  $S$ of $\D$ and a linear functional  $\boldsymbol{f}:\hbox{span}_\bbbr S\longrightarrow \bbbr$, the decomposition
\begin{equation}\label{tri-dec}
S=S^+\cup S^\circ\cup S^-
\end{equation}
where
 $$\hbox{\small $S^{\pm}=S^{\pm}_{\bs{f}}:=\{\a\in S\mid \boldsymbol{f}(\a)\gtrless0\}\andd S^\circ=S^\circ_{\bs{f}}:=\{\a\in S\mid \boldsymbol{f}(\a)=0\}$}$$
 is called a {\it triangular decomposition} for $S$.
It is called {\it trivial} if $S=S^\circ_{\bs{f}}.$
A subset $P$ of the root system ${{\D}}$ of $\fk$ is called a {\it parabolic} subset  of ${{\D}}$ if  $${{\D}}=P\cup -P\andd {{\D}}\cap (P+P)\sub P.$$

\begin{exa}\label{exa1}
{\rm Suppose that ${{\D}}={{\D}}^+\cup {{\D}}^\circ\cup {{\D}}^-$ is a nontrivial triangular decomposition  for ${{\D}}$ and  ${{\D}}^\circ={{\D}}^{\circ,+}\cup {{\D}}^{\circ,\circ}\cup {{\D}}^{\circ,-}$ is a  triangular decomposition for ${{\D}}^\circ$.  Then,  ${{\D}}^+\cup {{\D}}^{\circ,+}\cup {{\D}}^{\circ,\circ}$ is a parabolic subset of ${{\D}}.$
 }\end{exa}

\medskip

Suppose that $P$ is a parabolic subset of ${\D}$ and  set
\begin{equation}\label{p}
\displaystyle{\fk^\circ_{_P}:=\bigoplus_{\a\in P\cap-P}\fk^\a,\;\;\fk^+_{_P}:=\bigoplus_{\a\in P\setminus-P}\fk^\a,\;\; \fk^-_{_P}:=\bigoplus_{\a\in -P\setminus P}\fk^\a}\andd \fp:=\fk^\circ_{_P}\op\fk^+_{_P}.
\end{equation}

Assume  $\Omega$ is a module over $\fk^\circ_{_P}.$
 Consider $\Omega$ as a module over $\fp$ with trivial action of $\fk^+_P.$ Denoting  by $U(\cdot),$ the universal enveloping algebra of a Lie superalgebra, if the $\fk$-module $U(\fk)\ot_{U(\fp)}\Omega$ has a unique maximal submodule $Z$ intersecting $\Omega$ trivially, we set  $${\rm Ind}_P(\Omega):=\frac{U(\fk)\ot_{U(\fp)}\Omega}{Z}.$$ An irreducible finite weight  $\fk$-module $V$ is called {\it parabolically induced} if there is a parabolic subset $P$ of ${\D}$ and an irreducible module $\Omega$ over $\fk^\circ_{_P}$ such that $V\simeq {\rm Ind}_P(\Omega).$ An irreducible finite weight $\fk$-module $V$ is called {\it cuspidal} if it is not parabolically induced.

\section{Twisted Affine Lie superalgebras}\label{twisted affine}
Suppose that  $\fg$ is a  finite dimensional  basic classical simple Lie superalgebra of type $X=A(k,\ell) (\hbox{\tiny$(k,\ell)\neq (1,1)$}), D(k,\ell)$ with the standard Cartan subalgebra\footnote{It is the subalgebra of all diagonal matrices.} $H\sub\fg_0;$ here $k$ is a nonnegative integer and $\ell$ is a positive integer. Suppose that $\kappa\fm$ is  a nondegenerate supersymmetric invariant  even  bilinear form on $\fg.$ In \cite{van-thes}, the author introduces
a certain automorphism $\sg:\fg\longrightarrow \fg$  such that
\begin{itemize}
\item{\small \rm $\sg$ is of order $n=4$ if $X=A(2k,2\ell),$}
\item{\small \rm $\sg$ is of order $n=2$ if $X=A(2k-1,2\ell-1),A(2k,2\ell-1),D(k,\ell).$}
\end{itemize} Suppose $\zeta$ is the $n$-th primitive root of unity.
  Then,
we have
\[\fg=\bigoplus_{j=0}^{n-1}{}^{[j]}\fg\quad\quad\hbox{ where } {}^{[j]}\fg:=\{x\in\fg\mid \sg(x)=\zeta^jx\}\quad(0\leq j\leq n-1).\]
 For a two-dimensional vector space $\bbbc c\op\bbbc d,$ set
 \begin{equation}\label{affine}
 \fl:=\bigoplus_{j=0}^{n-1}({}^{[j]}\fg\ot t^j\bbbc[t^{\pm n}])\op\bbbc c\op\bbbc d\andd \fh:=(({}^{[0]}\fg\cap H)\ot 1)\op\bbbc c\op\bbbc d.
 \end{equation}
 Then  $\fl,$ which is denoted by $X^{(n)},$ together with
$$[x\ot t^p+rc+sd,y\ot t^q+r'c+s'd]:=[x,y]\ot t^{p+q}+p\kappa(x,y)\d_{p+q,0}c+sqy\ot t^q-s'px\ot t^p$$ (for $p,q\in\bbbz,$ $x,y\in \fg$ and $r,r',s,s'\in\bbbc$) in which ``~$\d_{\cdot,\cdot}$~" indicates the Kronecker delta,  is called  the {\it twisted affine Lie superalgebra} of type $X^{(n)}.$ We refer to the central element $c$ as the {\it canonical central element} of $\fl.$  For the details regarding affine Lie superalgebras see \cite{van-thes}; see also \cite[Appendix]{you8} and \cite{DSY}.

The Lie superalgebra $\fl$ has a  root space decomposition with respect to  $\fh.$ We denote  the corresponding root system by $R$ and  refer to $\fh$ as the {\it standard Cartan subalgebra} of $\fl.$ We have $R=R_0\cup R_1$ where $R_0$ (resp. $R_1$) is the set of weights of $\fl_0$ (resp. $\fl_1$) with respect to $\fh.$

The form $\kappa\fm$ induces the following nondegenerate  supersymmetric invariant  bilinear form  $\fm$  on $\fl:$
$$(x\ot t^p+rc+sd,y\ot t^q+r'c+s'd):=\kappa(x,y)\d_{p+q,0}+rs'+sr'.$$
 As the form is nondegenerate on  $\fh,$ one can transfer the  form  on  $\fh$ to a form on $\fh^*$ denoted again by $\fm.$
 The root system $R$ of $\fl$ with respect to $\fh$ is as in the following table:

 \begin{table}[h]\caption{} \label{table1}
 {\footnotesize \begin{tabular}{|c|l|}
\hline
$X^{(n)}$ &\hspace{3.25cm}$R$ \\
\hline
$A(2k,2\ell-1)^{(2)}$&$\begin{array}{rcl}
\bbbz\d
&\cup& \bbbz\d\pm\{\ep_i,\d_j,\ep_i\pm\ep_r,\d_j\pm\d_s,\ep_i\pm\d_j\mid i\neq r,j\neq s\}\\
&\cup& (2\bbbz+1)\d\pm\{2\ep_i\mid 1\leq i\leq k\}\\
&\cup& 2\bbbz\d\pm\{2\d_j\mid 1\leq j\leq \ell\}.
\end{array}$\\
\hline
$A(2k-1,2\ell-1)^{(2)},\;(k,\ell)\neq (1,1)$& $\begin{array}{rcl}
\bbbz\d&\cup& \bbbz\d\pm\{\ep_i\pm\ep_r,\d_j\pm\d_s,\d_j\pm\ep_i\mid i\neq r,j\neq s\}\\
&\cup& (2\bbbz+1)\d\pm\{2\ep_i\mid 1\leq i\leq k\}\\
&\cup& 2\bbbz\d\pm\{2\d_j\mid 1\leq j\leq \ell\}
\end{array}$\\
\hline
$A(2k,2\ell)^{(4)}$& $\begin{array}{rcl}
\bbbz\d&\cup&  \bbbz\d\pm\{\ep_i,\d_j\mid 1\leq i\leq k,\;1\leq j\leq \ell\}\\
&\cup& 2\bbbz\d\pm\{\ep_i\pm\ep_r,\d_j\pm\d_s,\d_j\pm\ep_i\mid i\neq r,j\neq s\}\\
&\cup&(4\bbbz+2)\d\pm\{2\ep_i\mid 1\leq i\leq k\}\\
&\cup& 4\bbbz\d\pm\{2\d_j\mid 1\leq j\leq \ell\}
\end{array}$\\
\hline
$D(k+1,\ell)^{(2)}$& $\begin{array}{rcl}
\bbbz\d&\cup&  \bbbz\d\pm\{\ep_i,\d_j\mid 1\leq i\leq k,\;1\leq j\leq \ell\}\\
&\cup& 2\bbbz\d\pm\{2\d_j,\ep_i\pm\ep_r,\d_j\pm\d_s,\d_j\pm\ep_i\mid i\neq r,j\neq s\}
\end{array}$\\
\hline
 \end{tabular}}
 \end{table}
\noindent in which  $\d,\ep_i,\d_p\in \sspan_\bbbc R$ $(1\leq i\leq k,\;1\leq p\leq \ell)$ satisfy
\[(\ep_i,\ep_j)=\d_{i,j},\;(\d_p,\d_q)=-\d_{p,q},\;(\ep_i,\d_p)=0\andd (\d,R)=\{0\}.\]
The root system $R_0$  of $\fl_0$ has a decomposition  $R_0=R(1)\cup R(2)$ with
{\small \begin{table}[h]\caption{}\label{table0}
{\footnotesize \begin{tabular}{|c|l|l|}
\hline
$X^{(n)}$& \multicolumn{1}{>{\centering\arraybackslash}m{45mm}|}{$R(1)$}  &\multicolumn{1}{>{\centering\arraybackslash}m{45mm}|}{$R(2)\hbox{ if $k\neq 0$}$}\\
\hline
$A(2k,2\ell-1)^{(2)}$
& $\begin{array}{l}
(2\d_{\ell,1}+(1-\d_{\ell,1}))\bbbz\d\\
\cup \bbbz\d\pm\{\d_j\pm\d_s\mid 1\leq j\neq s\leq \ell\}\\
\cup 2\bbbz\d\pm\{2\d_j\mid 1\leq j\leq \ell\}
\end{array}$
&$\begin{array}{l}
\bbbz\d\\
\cup \bbbz\d\pm\{\ep_i,\ep_i\pm\ep_r\mid 1\leq i\neq r\leq k\}\\
\cup (2\bbbz+1)\d\pm\{2\ep_i\mid 1\leq i\leq k\}
\end{array}$
\\
\hline
$\begin{array}{c}
A(2k-1,2\ell-1)^{(2)}\\(k,\ell)\neq (1,1)
\end{array}
$
& $\begin{array}{l}
(2\d_{\ell,1}+(1-\d_{\ell,1}))\bbbz\d \\
\cup \bbbz\d\pm\{\d_j\pm\d_s\mid 1\leq j\neq s\leq \ell\}\\
\cup 2\bbbz\d\pm\{2\d_j\mid 1\leq j\leq \ell\}
\end{array}$
& $\begin{array}{l}
(2\d_{k,1}+(1-\d_{k,1}))\bbbz\d\\
\cup \bbbz\d\pm\{\ep_i\pm\ep_r\mid 1\leq i\neq r\leq k\}\\
\cup (2\bbbz+1)\d\pm\{2\ep_i\mid 1\leq i\leq k\}
\end{array}$
\\
\hline
$A(2k,2\ell)^{(4)}$
& $\begin{array}{l}
2\bbbz\d\\
\cup (2\bbbz+1)\d\pm\{\d_j\mid 1\leq j\leq \ell\}\\
\cup 2\bbbz\d\pm\{\d_j\pm\d_s\mid 1\leq j\neq s\leq \ell\}\\
\cup 4\bbbz\d\pm\{2\d_j\mid 1\leq j\leq \ell\}
\end{array}$
& $\begin{array}{l}
2\bbbz\d\\
\cup  2\bbbz\d\pm\{\ep_i\mid 1\leq i\leq k\}\\
\cup 2\bbbz\d\pm\{\ep_i\pm\ep_r\mid 1\leq i\neq r\leq k\}\\
\cup(4\bbbz+2)\d\pm\{2\ep_i\mid 1\leq i\leq k\}
\end{array}$
\\
\hline
$
D(k+1,\ell)^{(2)}
$
& $\begin{array}{l}
2\bbbz\d\\
\cup 2\bbbz\d\pm\{\d_j\pm\d_s\mid 1\leq j, s\leq \ell\}
\end{array}$
& $\begin{array}{l}
\bbbz\d\\
\cup \bbbz\d\pm\{\ep_i\mid {1\leq i\leq k}\}\\
\cup 2\bbbz\d\pm\{\ep_i\pm\ep_r\mid 1\leq i\neq r\leq k\}
\end{array}$
\\
\hline
 \end{tabular}}
 \end{table}
}

 Also, we have  $R_{re}\cup R_{im}=S(1)\cup S(2)$ where $R(2)=S(2)=\emptyset$ if $k=0$,
 \begin{equation}\label{s(i)}
 S(1)=\bbbz\d\cup R(1)\cup (R\cap \frac{1}{2}R(1))\andd S(2)=\bbbz\d\cup R(2)\cup (R\cap \frac{1}{2}R(2)) \hbox{ if $k\neq 0$}.
 \end{equation}
 Recall that a subset $T\sub R$ is called {\it closed} if $(T+T)\cap R\sub T$ and set
 \begin{equation}\label{gi}
 \scg(i):=\bigoplus_{\a\in S(i)}\fl^\a\quad\quad(i=1,2).
\end{equation}
  Since $S(i)$ is a closed subset of $R$, we get that  $\scg(i)$ is a sub-superalgebra of $\fl.$

Using  (\ref{decom}), for each subset $T$ of $R,$ we set
\[T_*:=T\cap R_*\andd T^\times_*:=T\cap R^\times_*\quad\quad(*=re,im,ns), \]
\begin{equation}\label{set}\dot T:=\{\dot\a\in \sspan_\bbbc\{\ep_i,\d_p\mid1\leq i\leq k,~1\leq p\leq \ell\}\mid \exists \sg\in\bbbz\d\ni \dot\a+\sg\in T\}.\end{equation}
Also, we set
\begin{align}
&S_{\dot\a}:=\{\sg\in \bbbz\d\mid \dot\a+\sg\in R\}\quad\quad(\dot\a\in \dot R),\nonumber\\
&\dot T_{*}:=\{\dot\a\in \dot T\mid \dot\a+S_{\dot\a}\sub R_{*}\},\;\; \dot T_{*}^\times:=\dot T_*\setminus\{0\}\quad\quad(*=re,ns),\label{sets}\\
&\dot T_{sh}:=\{\dot\a\in \dot T_{re}^\times\mid |(\dot\a,\dot\a)|\leq |(\dot\b,\dot\b)|\;\;\;\forall\dot\b\in \dot T^\times_{re}\},\nonumber\\
&\dot T_{ex}:=2\dot T_{sh}\cap T\andd
\dot T_{lg}:=\dot T_{re}^\times\setminus(\dot T_{sh}\cup \dot T_{ex}).\nonumber
\end{align}
  By \cite[Table~4 and (3.9)]{you8}, we have:
\begin{align}
&(a)~\hbox{\footnotesize  $\forall\dot\a\in \dot{R}^\times$  $\exists~ r_{\dot\a}\in\{1,2,4\}$
 and $0\leq k_{\dot\a}<r_{\dot\a}$
 $\ni S_{\dot\a}=(r_{\dot\a}\bbbz+k_{\dot\a})\d,$
 }\nonumber\\
&(b)~ \hbox{\small  $\forall \dot\a,\dot\b\in \dot R_{sh},$
$S_{\dot\a}=S_{\dot\b}=\bbbz\d$ and $\forall\dot\gamma,\dot\eta\in \dot R_{ns}^\times$   $\exists ~s\in\bbbz^{>0}~\ni S_{\dot\gamma}=S_{\dot\eta}=s\bbbz\d$},\label{imp}\\
&(c)~\hbox{\footnotesize
 {if $T$ is a  closed subset  of $R$,  $\bbbz\d\sub T$ and  $\dot\a\in \dot T,$ we have
$\dot\a+S_{\dot\a}\sub T$.}}\nonumber
\end{align}
\section{quasi-Integrable Modules}
%
%
%
%
%

\begin{deft}{\rm Using the same notations as in \S~\ref{twisted affine}, we suppose that $\fl$ is a twisted affine Lie superalgebra with $R(2)\neq \emptyset.$
\begin{itemize}
\item[(i)] For an irreducible finite weight module $M$ over $\fl,$ we say a  subset $S$ of $R$ is  {\it tight} if there is a  nonzero real root $\a\in S$ with $(\a+\bbbz\d)\cap S\sub S^{ln}(M)$ or   $(\a+\bbbz\d)\cap S\sub  S^{in}(M);$ otherwise, we call it {\it hybrid}.
\item[(ii)] Recall $\scg(i)$ ($i=1,2$) from (\ref{gi}). Suppose $M$ is an irreducible  $\fl$-module. The $\fl$-module  $M$ is called  {\it integrable} if $R_{re}^\times=R^{ln}.$ If $\{r,t\}=\{1,2\},$ the irreducible $\fl$-module $M$ is called {\it $t$-quasi-integrable} if $S(r)$ is  hybrid and $M$  is integrable as a module over  $\scg(t),$  that is $S(t)\cap R_{re}^\times\sub R^{ln}$.
    \end{itemize}}
\end{deft}

\begin{exa}
{\rm
 In this example, we present a module over $\fl:=A(2k-1,1)^{(2)}$  $(k\geq2)$ which is $2$-quasi integrable but it is not a highest weight module. Suppose $\fh$ is  the standard Cartan subalgebra of $\fl$ with corresponding root system $R.$

We make a convention that in what follows  for a subalgebra $\ft$ of $\fl$ containing $\fh$ with  corresponding root system $\frak{s}$, a  parabolic subset $P$ of $\frak{s}$ and an irreducible module $\Omega$ over $\ft^\circ_{_P},$ the module action on ${\rm Ind}_P(\Omega)$ is denoted, as usual, by juxtaposition while the module action on $U(\ft)\ot_{_{U(\ft^\circ_{_P})}}\Omega$ is denoted by $``*".$

We recall that the root system $R$ of $A(2k-1,1)^{(2)}$ ($k\geq 2$) is
\begin{align*}
R=\bbbz\d&\cup (\pm 2\d_1+2\bbbz\d)\cup  (\{\pm\ep_i\pm\d_1\mid 1\leq i\leq k\}+\bbbz\d)\\
&\cup (\{\pm 2\ep_i\mid 1\leq i\leq k\}+2\bbbz\d+\d)  \cup (\{\pm\ep_i\pm\ep_j\mid 1\leq i<j\leq k\}+\bbbz\d).
\end{align*}
We define
\begin{align*}
{\small P_1}:=&{\small\{0,\ep_k\pm\d_1,\pm2\d_1\}},\\
{\small P_2}:=&{\small\{0\}\cup \{\ep_i\pm\ep_j\mid 1\leq i<j\leq k\}\cup \{\ep_i\pm\d_1\mid 1\leq i\leq k-1\}\cup \{\pm\ep_k\pm\d_1\}\cup \{\pm2\d_1\},}\\
{\small P_3}:=& {\small \{0\}\cup (\{\pm\ep_i\pm\ep_j\mid 1\leq i\neq j\leq k\}+\bbbz^{\geq 0}\d)\cup  (\{\pm2\ep_i\mid 1\leq i\leq k\}+\d+2\bbbz^{\geq 0}\d)}\\
&{\small\cup (\{\pm2\d_1\}+2\bbbz^{\geq 0}\d)\cup (\{\pm\ep_i\pm\d_1\mid 1\leq i\leq k\}+\bbbz^{\geq 0}\d)}.
\end{align*}
We have
\begin{align*}
\frak{s}_1:=&P_1\cap -P_1=\{0,\pm2\d_1\}~~\hbox{ (root system of $A_1$)},\\
\frak{s}_2:=&P_2\cap -P_2=\{0,\pm\ep_k\pm\d_1,\pm2\d_1\}~~\hbox{ (root system of $C(2)$)},\\
\frak{s}_3:=&P_3\cap -P_3=\{0,\pm2\d_1,\pm\ep_i\pm\ep_j,\pm\ep_i\pm\d_1\mid 1\leq i\neq j\leq k\}~~\hbox{ (root system of $D(k,1)).$}
\end{align*} Note that $P_i$ is a parabolic subset of $\frak{s}_{i+1}$ ($i=1,2$) and $P_3$ is a parabolic subset of $R.$
Recall (\ref{talpha}) as well as (\ref{affine}) and set
\[\hB_i:=\bbbc c+\bbbc d+\bbbc t_{2\d_1}\op\sum_{j=1}^k\bbbc t_{\ep_j}\op\bigoplus_{0\neq \a\in \frak{s}_i}\fl^\a\andd\hB_i^+=\bigoplus_{0\neq \a\in P_{i-1}\setminus -P_{i-1}}\fl^\a\quad\quad(i=2,3).\]

Assume $\zeta\in\bbbq\setminus\bbbz$ and $\xi$ is an irrational number. Fix $e\in\fl^{2\d_1}$ and $f\in\fl^{-2\d_1}$ with $[e,f]=t_{2\d_1};$ see  (\ref{talpha}) and \cite[Lem.~3.1]{you3}. Suppose  $K_1$ is the vector space with  basis $\{v_\mu\mid \mu\in \zeta+2\bbbz\}.$ Then $K_1$ together with

$$\begin{array}{rlrlll}
\hspace{-1mm}d v_\mu=&\hspace{-3mm}0,& \hspace{-1mm}c v_\mu=&\hspace{-3mm}(2k+2)v_\mu, & f v_\mu=\frac{-1}{2}(\xi-(\mu-1)^2)v_{\mu-2},&\hspace{-1mm}e v_\mu= v_{\mu+2},\\
\hspace{-1mm}t_{2\d_1} v_\mu=& \hspace{-3mm}-2\mu v_\mu,&\hspace{-1mm} t_{\ep_i}v_\mu=&\hspace{-3mm}(k-i+2)v_\mu &(1\leq i\leq k,\mu\in \zeta+2\bbbz)&
\end{array}$$
 is an irreducible finite weight   module over the finite dimensional reductive  Lie algebra $\hB_1:=\bbbc c\op\bbbc d\op\bbbc t_{2\d_1}\op\sum_{i=1}^k\bbbc t_{\ep_i}\op\bbbc e\op\bbbc f$ with respect to $\bbbc c\op\bbbc d\op\bbbc t_{2\d_1}\op\sum_{i=1}^k\bbbc t_{\ep_i}.$ We note that both $e$ and $f$ act injectively on $K_1,$ so
\[\pm2\d_1\in \frak{s}_1^{in}(K_1).\] Define
\begin{align*}
\varrho:&\bbbc c\op\bbbc d\op\bbbc t_{2\d_1}\op \sum_{i=1}^k \bbbc t_{\ep_i}\longrightarrow \bbbc\\
& c\mapsto 2k+2,\; d\mapsto0,\; t_{2\d_1}\mapsto -2\zeta,\; t_{\ep_{i}}\mapsto k-i+2\;\;\;(1\leq i\leq k).
\end{align*}Then
\begin{equation}
\label{supp}
 \supp(K_1)=\varrho+2\bbbz\d_1.
 \end{equation}
We next set
\[K_i:={\rm Ind}_{_{P_{i-1}}}(K_{i-1})\quad\quad(i=2,3,4)\] and note that
\[K_1\lhook\joinrel\xrightarrow{\hbox{\tiny as $\hB_1$-submod}} K_2\lhook\joinrel\xrightarrow{\hbox{\tiny as $\hB_2$-submod}} K_3\lhook\joinrel\xrightarrow{\hbox{\tiny as $\hB_3$-submod}} K_4.\] Since $K_1$ is an irreducible finite weight module over $\hB_1,$ it is not hard to see that  $K_i$ ($i=2,3$) is also an irreducible finite weight  module over $\hB_i$ and $M:=K_4$ is an irreducible finite weight module over $\fl;$ see \cite[Pro.~1.8]{DG} and \cite[Lem.~2.3]{DMP}. The aim of this example is showing that $M$ is $2$-quasi-integrable while it is not a highest weight module.

For  a nonzero nonsingular root $\a,$ we have $2\a\not\in R$ and $\dim(\fl^\a)=1$. So,  we get
\begin{equation}\label{zero-1}\fl^\a\fl^\a v=[\fl^\a,\fl^\a]v=\{0\}\quad \hbox{ for } \a\in R_{ns}^\times\cap \frak{s}_i,\;v\in K_i~(i=1,2,3).\end{equation}
Also from the structure of $K_i$ ($i=1,2,3$), we have
\begin{equation}\label{zero}\fl^\a* K_i=\{0\}\quad\quad \hbox{ for }\a\in P_i\setminus-P_i~(i=1,2,3).\end{equation}

\noindent{\bf Step~1.} $\supp(K_2)\sub \varrho+\bbbz2\d_1-\{0,1,2\}\ep_k:$  Since $K_2={\rm Ind}_{_{P_1}}(K_1)$ is generated by $K_1$,  we are done  using (\ref{supp}), (\ref{zero-1}), (\ref{zero}) together with the fact that $-P_1\setminus P_1=\{-\ep_k\pm\d_1\}\sub (\frak{s}_2)_{ns}^\times$.

\noindent{\bf Step~2.}  $\{\pm(\ep_1\pm\ep_{2}),\ldots,\pm(\ep_{k-1}\pm\ep_k)\}\sub \frak{s}_3^{ln}(K_3):$ One knows that
\[B:=\{\a_j:=\ep_j-\ep_{j+1},\a_k:=\ep_k-\d_1,\a_{k+1}:=2\d_1\mid 1\leq j\leq k-1\}\]is  a base of $\frak{s}_3=D(k,1).$ For each $1\leq j\leq k+1,$ fix $ e_{\a_j}\in\fl^{\a_j}$ and $f_{\a_j}\in\fl^{-\a_j}$ with $[e_{\a_j},f_{\a_j}]=t_{\a_j}.$ For $1\leq j\leq k-1$ and a weight vector  $v\in  K_1\sub K_2$ of weight $\lam\in \varrho+2\bbbz\d_1$, as $B$ is a base, using PBW-Theorem, the $\hB_3$-submodule of $U(\hB_3)\ot_{U(\hB_2\op\hB_3^+)}K_2$ generated by $f_{\a_j}^2*v$
is the linear span of the elements of the form
\[f_{\a_{p_1}}*\cdots *f_{\a_{p_s}}*e_{\a_{j_1}}*\cdots* e_{\a_{j_t}}*f_{\a_j}^2 *v,\; f_{\a_{p_1}}*\cdots *f_{\a_{p_s}}*f_{\a_j}^2*v,\;e_{\a_{j_1}}*\cdots * e_{\a_{j_t}}*f_{\a_j}^2* v,\;f_{\a_j}^2*v\] for $1\leq j_1,\ldots,j_t,p_1,\ldots,p_s\leq k+1.$
Since $j\in\{1,\ldots, k-1\}$ and  for $p\neq k+1,$ by (\ref{zero}), $e_p*v=0,$   we have for $1\leq p\leq k+1$ that
\[\hbox{\footnotesize $e_{\a_p}f^2_{\a_j}*v=
\left\{
\begin{array}{ll}
f_{\a_j}^2e_p*v=0& p\neq j,k+1,\\
(2\varrho-\a_j)(t_{\a_j})f_{\a_j}*v=0& p=j,\\
f_{\a_j}^2e_{2\d_1}*v\in f_{\a_j}^2*K_1^{\lam+2\d_1}&p=k+1,
\end{array}
\right.$}
\]
that is
\[e_{\a_p}*f^2_{\a_j}*K_1^\lam=e_{\a_p}f^2_{\a_j}*K_1^\lam\sub f^2_{\a_j} *K_1^{\lam+2\d_1}.\] So, we get that  \[e_{\a_{j_1}}*\cdots * e_{\a_{j_t}}*f_{\a_j}^2* K_1\sub f_{\a_j}^2* K_1.\] This implies that
the support of the  $\hB_3$-module generated by $f^2_{\a_j}*v$  ($1\leq j\leq k-1$) lies in $\varrho+\bbbz\d_1-2(\ep_j-\ep_{j+1})+\sum_{t=1}^{k}\bbbz^{\leq 0}\a_t.$ This, in particular, together with Step~1 implies that the $\hB_3$-submodule of $U(\hB_3)\ot_{U(\hB_2\op\hB_3^+)}K_2$ generated by $f^2_{\a_j}*v$ intersects $K_2$ trivially and so $f^2_{\a_j}v=0$ $(1\leq j\leq k-1);$ in other words,
\begin{equation}\label{alphj}
-\a_j=-(\ep_j-\ep_{j+1})\in \frak{s}_3^{ln}(K_3)\quad\quad(1\leq j\leq k-1).
\end{equation}
Next, note that
\[B':=\{\b_i:=\ep_i-\ep_{i+1},\b_{k-1}:=\ep_{k-1}+\ep_k,\b_k:=-\ep_k-\d_1,\b_{k+1}:=2\d_1\mid 1\leq i\leq k-2\}\] is also a base  of $\frak{s}_3$.
Contemplating Step~1 and suppose \[r:={\rm max}\{j\mid \exists n\in\bbbz\ni \varrho-j\ep_k+n\d_1\in \supp(K_2)\}\andd \mu:=\varrho-r\ep_k+m\d_1\in\supp(K_2)\] for some $m\in\bbbz.$ Fix $0\neq v\in K_2^\mu=K_2^{\varrho-r\ep_k+m\d_1}.$ For $1\leq i\leq k+1,$ fix $ e_{\b_i}\in\fl^{\b_i}$ and $f_{\b_i}\in\fl^{-\b_i}$ with $[e_{\b_i},f_{\b_i}]=t_{\b_i}.$
 Setting
\[s:=\mu(t_{\b_{k-1}})=(\varrho-r\ep_k+m\d_1)(t_{\b_{k-1}})=\varrho(t_{\ep_{k-1}})+\varrho(t_{\ep_k})-r=5-r\geq0,\] as $e_{\b_k}*v\in K_2^{\mu-\ep_k-\d_1}=K_2^{\varrho-(r+1)\ep_k+(m-1)\d_1}=\{0\}$,
we have using (\ref{zero}) that
\[e_{\b_t}f^{s+1}_{\b_{k-1}}*v=\left\{
\begin{array}{ll}
f^{s+1}_{\b_{k-1}}e_{\b_t}*v=0&t=k,~1\leq  t< k-2,\\
(s+1)(\mu(t_{\b_{k-1}})-s)f_{\b_{k-1}}^sv=0& t=k-1,\\
f^{s+1}_{\b_{k-1}}e_{2\d_1}*v& t=k+1.
\end{array}
\right.
\]
 So for $1\leq j_1,\ldots,j_t\leq k,$ we get \[e_{\b_{j_1}}\cdots e_{\b_{j_t}}f^{s+1}_{\b_{k-1}} *v\in f^{s+1}_{\b_{k-1}} *\sum_{j\in\bbbz}K_2^{\varrho+j\d_1-r\ep_k}.\]
This implies that
the support of the  $\hB_3$-submodule of $U(\hB_3)\ot_{U(\hB_2\op\hB_3^+)}K_2$ generated by
$f^{s+1}_{\b_{k-1}}*v$ lies in $\varrho+\bbbz\d_1-r\ep_k-(s+1)\b_{k-1}+\sum_{t=1}^{k}\bbbz^{\leq 0}\b_t.$ This together with Step~1 implies that the $\hB_3$-submodule generated by $f^{s+1}_{\b_{k-1}}*v$ intersects $K_2$ trivially and so $f^{s+1}_{\b_{k-1}}v=0.$ This means that $-\b_{k-1}=-\ep_{k-1}-\ep_{k}\in \frak{s}_3^{ln}(K_3).$ Since $\{\ep_i-\ep_{i+1},\ep_{k-1}+\ep_k\mid 1\leq i\leq k-1\}$ is a base of $D(k),$ this together with (\ref{alphj}), (\ref{zero}) and \cite[Thm.~4.7]{you8} completes the proof of this step.

\noindent{\bf Step~3.}  $\pm\{2\d-(\ep_1+\ep_2),\ep_i-\ep_{i+1},2\ep_k-\d\mid 1\leq i\leq k-1\}\sub R^{ln}(K_4):$ We first mention that $\hbox{\small$\D:=\{\gamma_1:=-2\d_1,\gamma_2:=\d_1+\ep_k,\gamma_{i+1}:=-\ep_{i}+\ep_{i-1},
\gamma_{k+2}:=\d-2\ep_1\mid 2\leq i\leq k\}$}$ is a linearly independent subset of $R$ with $R\sub \sspan_{\bbbz^{\geq 0}}\D\cup \sspan_{\bbbz^{\leq 0}}\D.$
 Suppose that $0\neq v\in K_1\sub K_2\sub K_3$ is a weight vector of weight $\lam\in \varrho+2\bbbz\d_1.$ We note that as
\[\D\setminus\{-2\d_1\}\sub (P_3\setminus-P_3)\cup (P_2\setminus-P_2)\cup (P_1\setminus-P_1),\] we have using (\ref{zero}) that
$e_{\gamma_i}*v=0$ $(2\leq i\leq k+2).$ So,
 for $1\leq p\leq k+1,$ we have
\[e_{\gamma_p}f_{\gamma_{k+2}}*v=\left\{
\begin{array}{ll}
f_{\gamma_{k+2}}e_{\gamma_p}*v=0&2\leq  p\leq k+1,\\
\lam(t_{\gamma_{k+2}})v=0& p=k+2,\\
f_{\gamma_{k+2}}\underbrace{e_{-2\d_1}*v}_{\in K_1}& p=1.
\end{array}
\right.
\]
Setting  $\fl^+=\op_{\a\in P_3\setminus-P_3}\fl^\a$ and  using the same argument as before, one can get that the $\fl$-submodule of $U(\fl)\ot_{U(\hB_3\op\fl^+)}K_3$ generated by $f_{\gamma_{k+2}}*v$ intersects $K_3$ trivially and so $f_{\gamma_{k+2}}v=0,$ that is $-(\d-2\ep_1)=-\gamma_{k+2}\in R^{ln}(K_4).$
This together with Step~2, (\ref{zero}), \cite[Thm.~4.7]{you8} and the facts that
\begin{align*}
&\d-2\ep_k=2(\ep_{k-1}-\ep_k)+\cdots+2(\ep_1-\ep_2)+(\d-2\ep_1)\andd\\
& 2\d-(\ep_1+\ep_2)=(\d-2\ep_1)+(\ep_1+\ep_2)+2(\ep_1-\ep_2)+(\d-2\ep_1)
 \end{align*}
implies that $\pm(\d-2\ep_k),\pm(2\d-(\ep_1+\ep_2))\in R^{ln}(K_4).$ This completes the proof of this step due to Step~2.

\noindent{\bf Step~4.} $M=K_4$ is a $2$-quasi-integrable $\fl$-module which is not a highest weight module: Since by Step~3,
$\pm\{2\d-(\ep_1+\ep_2),\ep_i-\ep_{i+1},2\ep_k-\d\mid 1\leq i\leq k-1\}\sub R^{ln}(M),$ using \cite[Thm.~4.7]{you8}, one can easily see that for  $$\a\in R(2)_{re}^\times=(\{\pm\ep_i\pm\ep_r\mid 1\leq i\neq r\leq k\}+\bbbz\d)\cup(\{\pm2\ep_i\mid 1\leq i\leq k\}+\d+2\bbbz\d),$$ $\pm\a\in R^{ln},$ that is $R(2)_{re}^\times\sub R^{ln}(M)$ which in turn implies that $S(2)_{re}^\times\sub R^{ln}(M).$ Since $\{\pm2\d_1\}+\bbbz^{>0}\d\sub  P_3\setminus-P_3,$ we have $\{\pm2\d_1\}+\bbbz^{>0}\d\sub R^{ln}(M),$  while $\pm2\d_1\in R^{in}(M).$ In other words, $S(1)$
is hybrid. Therefore $M$ is a $2$-quasi integrable $\fl$-module. On the other hand, as both $2\d_1$ and its opposite $-2\d_1$ belong to $R^{in}(M),$ $M$ is not a highest weight module.}
\end{exa}

From now on till  end of the paper, we assume
{\it
$\fl$ is a twisted affine Lie superalgebra of one of the types $A(2k-1,2\ell-1)^{(2)}$  ({\tiny$(k,\ell)\neq (1,1)$}), $A(2k,2\ell)^{(4)},$  $A(2k,2\ell-1)^{(2)}$ or $D(k,\ell)^{(2)},$ where $k,\ell$ are ``positive integers". We keep the same notations as in \S~\ref{twisted affine}; in particular, we frequently use Tables~\ref{table1},\ref{table0} and (\ref{s(i)}) without further reference.}
We assume $M$ is an irreducible  finite weight module over $\fl.$
 Using \cite[Pro.~4.4]{you8}, $M$ has shadow. Moreover,  the canonical central element $c$ acts on $M$ as an scalar. This scalar is called the {\it level} of $M$.

\begin{rem}
{\rm One knows that if $\a,2\a\in R_{re}^\times$ and $*=ln,in,$ then $\a\in R^{*}$ if and only if $2\a\in R^{*};$ see \cite[Thm.~4.7]{you8} and note that $M$ has shadow. So if $j\in\{1,2\}$ and $S(j)$ is hybrid, then by \cite[\S~5]{you8}, there is $r\in\{\pm 1\}$ such that
 \[\hbox{\small $R(j)_{re}^\times=\{\a\in R(j)_{re}^\times\mid \exists N\in\bbbz^{>0}\;\ni \; (\a+\bbbz^{\geq N}r\d)\cap R(j)\sub  R^{ln}\}.$}\]
This together with the  fact that $S(j)_{re}\setminus R(j)\sub R_{re}\cap \frac{1}{2}R_{re},$ implies that
 \[\hbox{\small $S(j)_{re}^\times=\{\a\in S(j)_{re}^\times\mid \exists N\in\bbbz^{>0}\;\ni \;(\a+\bbbz^{\geq N}r\d)\cap S(j)\sub  R^{ln}\}$}.\]
If $r=1$ (resp. $r=-1$), we call $S(j)$ {\it up-nilpotent hybrid} (resp.  {\it down-nilpotent hybrid}).}

\end{rem}

Contemplating \cite[Lem.~3.4]{you8} as well as Lemma~2.4, Corollary~2.5 and Theorem~2.6 of \cite{you3} and using the super version of the standard techniques in Lie theory, we get the following proposition:
\begin{Pro}\label{gen-inf}
Suppose that  $T$ is a closed subset of $R.$ Set $\scg:=\op_{\a\in T}\fl^\a.$ Assume  $W$ is a $\scg$-submodule of $M$ and  $\a,-\a\in T\cap R^{ln}.$ Then we have the following:
\begin{itemize}
\item[(i)] If $\lam,\lam+\a\in \supp(W),$ then $\fl^\a W^\lam\neq \{0\}.$
\item[(ii)] For $\lam\in \supp(W)$, $\frac{2(\lam,\a)}{(\a,\a)}\in\bbbz.$ Moreover, if $\lam\in\supp(W)$ and  $\frac{2(\lam,\a)}{(\a,\a)}\in \bbbz^{>0}$ (resp. $\in \bbbz^{<0}$), then $\lam-\a\in \supp(W)$ (resp. $\lam+\a\in \supp(W)$).
\end{itemize}
\end{Pro}
\begin{Pro}
\label{hybrid}
 Let $T$ be a closed subset of $R\setminus R_{ns}^\times$ with $\bbbz\d\sub T$ and set $\scg:=\op_{\a\in T}\fl^\a$ (in particular, $\fh=\fl^0\sub \scg$). Suppose $W$ is a $\scg$-submodule of $M.$
\begin{itemize}
\item[(i)] Assume  $\aa\sub \supp(W)$ is nonempty and $S\sub T$ is a nonempty finite subset with
$$S\cap R_{im}=\emptyset,~S\sub \frak{B}_W,~ -S\sub \frak{C}_W\andd (\aA+S)\cap \supp(W)\sub \aA.$$ Then,  there is $\lam\in \aA$ such that $(\lam+\sspan_{\mathbb{Z}^{\geq 0}}S)\cap \supp(W)=\{\lam\}.$
\item[(ii)] Suppose that $\boldsymbol{f}$  is a functional on $\sspan_\bbbr T$ with $\boldsymbol{f}(\d)\neq 0$ and corresponding triangular decomposition $T=T^+\cup T^\circ\cup T^-.$ Assume for $K:=T_{re}\setminus \frac{1}{2}T_{re},$   $J:=K\cup (\bbbz\d\cap \sspan_\bbbz K)$ is an affine root system. If  $$T^+_{re}:=T^+\cap R_{re}\sub \frak{B}_W\andd T^-_{re}:=T^-\cap R_{re}\sub \frak{C}_W,$$  we have the following:
\begin{itemize}
\item[(a)]     Let $r\in\{\pm1\}$ with $\boldsymbol{f}(r\d)>0,$ then there are $\lam\in \supp(W)$ and  a positive integer $p$   with $(\lam+\bbbz^{>0}p(r\d))\cap\supp(W)=\emptyset,$
\item[(b)]  Recall (\ref{sets}) and (\ref{imp}(c)). If there is $\dot \a\in \dot T_{re}^\times$ with $S_{\dot\a}=\bbbz\d,$ then  there exists  $\mu\in\supp(W)$ such that $\mu+\a\not\in \supp(W)$ for all $\a\in  T^+.$
\end{itemize}
\end{itemize}
\end{Pro}
\pf (i) See \cite[Pro.~3.6]{you8}.

(ii)(a) follows from (i) together with the same argument as in \cite[Lem.~5.1]{you8}.

(ii)(b) follows from the same argument as in \cite[Pro.~3.7]{you8}\footnote{In Proposition~3.7 of \cite{you8}, we are working with $\fl$ or $\fl_0$ instead of $\scg$ introduced here, but our assumptions on $T$ make the same situation as in Proposition~3.7 of \cite{you8}.} by using part (ii)(a) and changing  the role of $\d$ with $-\d$ if $\boldsymbol{f}(\d)<0$.
\qed

\begin{Pro}\label{submod}
Recall   (\ref{gi}) and suppose $M$ is of nonzero level. Suppose  $W\sub M$ and $\{i,j\}= \{1,2\}.$
\begin{itemize}
\item[(i)] Let  $\boldsymbol{f}$ be a linear functional on $\sspan_\bbbr S(i)$ with corresponding  triangular decomposition
$S(i)=S(i)^+\cup S(i)^\circ\cup S(i)^-.$ If $W$ is a $\scg(i)$-submodule of $M$ and   $$\bs{f}(\d)\neq 0,~~S(i)^+\cap R_{re}\sub R^{ln},~~ S(i)^-\cap R_{re}\sub R^{in},$$ then there is $\mu\in\supp(W)$ such that $\mu+\a\not\in \supp(W)$ for all $\a\in  S(i)^+.$\\
\item[(ii)]   Assume $W$ is an integrable $\scg(j)$-submodule of $M$ and $\lam\in \supp(W)$ satisfies $\supp(W)\cap (\lam+\bbbz^{>0}r\d)=\emptyset$ for some $r\in\{\pm1\}.$ Then for each $\dot\a\in \dot S(j)$, see (\ref{set}), there is $m_{\dot\a}\in\bbbz^{\geq 0}$ such that $\lam+\dot\a+nr\d\not\in \supp(W)$ for all $n\geq m_{\dot\a}.$\\
\item[(iii)] Suppose that  $S(i)$ is up-nilpotent  (resp. down-nilpotent) hybrid. If  $W$ is a invariant under the action of  both $\scg(i)$ and  $\scg(j)$ such that it is integrable as a module over $\scg(j),$ then
 there are $\mu\in\supp(W)$ and a functional $\boldsymbol{f}$ on $\sspan_\bbbr R$  with  $\boldsymbol{f}(\d)>0 $ (resp. $\boldsymbol{f}(\d)<0 $) and  $\boldsymbol{f}(\dot S(j))=\{0\}$ such that
    \begin{itemize}
\item[$\bullet$] $\mu+\a\not\in \supp(W)$ for all $\a\in  S(i)$ with $\boldsymbol{f}(\a)>0,$
\item[$\bullet$] for each $\dot\a\in \dot S(j),$ there is $m_{\dot\a}\in\bbbz^{>0}$ such that $\mu+\dot\a+n\d\not\in \supp(W)$ (resp. $\mu+\dot\a-n\d\not\in \supp(W)$) for all $n\geq m_{\dot\a}.$
    \end{itemize}
\end{itemize}
\end{Pro}
\pf (i) We note that $S(i)^{in}\sub \frak{C}_W.$ Since   $M$ has  shadow, $R^{ln}\sub \frak{B}_M \sub \frak{B}_W.$ Therefore, the result follows from (\ref{s(i)})  and  Proposition~\ref{hybrid}(ii)(b).

(ii) We first assume $r=1.$
Suppose $\lam$ is as in the statement and $\dot\a\in  \dot S(i)^\times.$ Recalling (\ref{imp}), for each  $n\in\bbbz,$ we have
\begin{equation}\label{new03}
\dot\a+(k_{\dot\a}+r_{\dot\a}n)\d\in R_{re}^\times.
\end{equation} This in particular together with Proposition~\ref{gen-inf}(iii) implies that  \begin{equation}\label{new02}
{2(\dot\a+k_{\dot\a}\d,\lam)}/{(\dot\a,\dot\a)},{2(\dot\a+(k_{\dot\a}+r_{\dot\a})\d,\lam)}/{(\dot\a,\dot\a)}\in\bbbz
\end{equation} which in turn implies that ${2(r_{\dot\a}\d,\lam)}/{(\dot\a,\dot\a)}\in\bbbz.$
Since the canonical  central element  $c$ acts on $M$ as a nonzero scalar, we get $$(\lam,\d)M^\lam=\lam(c)M^\lam=cM^\lam\neq \{0\}.$$ Therefore,  $(\lam,\d)\neq 0$ and so
${2(r_{\dot\a}\d,\lam)}/{(\dot\a,\dot\a)}\in\bbbz\setminus\{0\}.$
\smallskip

 \noindent{{\bf Case~1.} ${2(\lam,r_{\dot\a}\d)}/{(\dot\a,\dot\a)}\in \bbbz^{>0}:$}
For each $\dot \a\in \dot S(i)^\times,$ considering (\ref{new02}), we choose $N\in \bbbz^{>0}$ such that ${2(\lam,\dot\a+(k_{\dot\a}+r_{\dot\a}N)\d)}/{(\dot\a,\dot\a)}$ is positive. We claim that for $n>k_{\dot\a}+r_{\dot\a}N,$ $\lam+\dot\a+n\d\not \in\supp(W).$ In fact, if $n>k_{\dot\a}+r_{\dot\a}N$ and $\lam+\dot\a+n\d\in\supp(W),$ since \[
{2(\lam+\dot\a+n\d,\dot\a+(k_{\dot\a}+r_{\dot\a}N)\d)}/{(\dot\a,\dot\a)}=({2(\lam,\dot\a+(k_{\dot\a}+r_{\dot\a}N)\d)}/{(\dot\a,\dot\a)})+2>0,\]
we get from Proposition~\ref{gen-inf}(iii) that $\lam+({n-k_{\dot\a}-r_{\dot\a}N})\d\in\supp(W)$ which is a contradiction.

\smallskip

 \noindent{{\bf Case~2.} ${2(\lam,r_{\dot\a}\d)}/{(\dot\a,\dot\a)}\in \bbbz^{<0}:$}
We claim that $\lam+\dot\a+n\d\not \in \supp(W)$ for all  $\dot\a\in \dot S(i)^\times$ and all  positive integers $n.$ To the contrary, assume $\lam+\dot\a+n\d\in \supp(W)$ for some  $\dot\a\in \dot S(i)^\times$ and some positive integer $n.$ Contemplate (\ref{new02}) and (\ref{new03}) and choose $p\in\bbbz^{>0}$ such that $pr_{\dot\a}>k_{\dot\a}$ and   ${2(\lam,\dot\a+r_{\dot\a}(-p)\d+k_{\dot\a}\d)}/{(\dot\a,\dot \a)}>0.$ This gives that
\[pr_{\dot\a}-k_{\dot\a}>0\andd {2(\lam+\dot\a+n\d,\dot\a+r_{\dot\a}(-p)\d+k_{\dot\a}\d)}/{(\dot\a,\dot\a)}>0,\] and so  Proposition~\ref{gen-inf}(iii) implies that
$$\lam+({n+pr_{\dot\a}-k_{\dot\a}})\d=\lam+\dot\a+n\d-(\dot\a+r_{\dot\a}(-p)\d+k_{\dot\a}\d)\in \supp(W)$$ which is a contradiction.

A mild modification of what we did for $r=1$ gives  the result in case $r=-1.$

\smallskip

(iii)
Set  $P:=S(i)^{ln}\cup -S(i)^{in}\cup \bbbz^{\geq 0}\d$ if $S(i)$ is up-nilpotent hybrid and  $P:=S(i)^{ln}\cup -S(i)^{in}\cup \bbbz^{\leq 0}\d$ if $S(i)$ is down-nilpotent hybrid. Then $P$  is a parabolic subset of $S(i);$ see the proof of \cite[Lem.~5.4]{you8}.
Since $R(i)=(S(i)_{re}\setminus \frac{1}{2}S(i)_{re})\cup (\bbbz\d\cap \sspan_\bbbz S(i)_{re})$ is the root system of an affine Lie algebra and $R(i)\cap P$ is a parabolic subset of $R(i)$,  by \cite[Pro.~2.10(ii)]{DFG}, there is a functional $\boldsymbol{f}$ on $\sspan_\bbbr R(i)$ such that $$R(i)\cap P=\{\a\in R(i)\mid \boldsymbol{f}(\a)\geq 0\}.$$ This in turn implies that
$P=\{\a\in S(i)\mid \boldsymbol{f}(\a)\geq 0\}.$
Extend $\boldsymbol{f}$ on $\sspan_\bbbr R$ with $\boldsymbol{f}(\dot S(j))=\{0\},$ we get the result using parts~(i),(ii).
\qed
\begin{Pro}
\label{gen}
Suppose that $\{i,j\}=\{1,2\}$ and assume $M$ has nonzero level.
Assume $S(i)$ is up-nilpotent  (resp. down-nilpotent) hybrid and $M$ is  integrable as a  module over $\scg(j)$.   Then, there are a nonzero weight vector  $v\in M$  and  a functional $\boldsymbol{f}$ on $\sspan_\bbbr R$ with corresponding triangular decomposition $R=R^+\cup R^\circ\cup R^-$ such that
\begin{itemize}
\item[$\bullet$] $\boldsymbol{f}(\d)>0$ (\hbox{resp.} $\bs{f}(\d)<0$)
and $\boldsymbol{f}(\dot S(j))=\{0\},$
\item[$\bullet$] $\forall~  \a\in R^+\cap S(i),\;\;  \fl^{\a} v=\{0\},$
\item[$\bullet$]
$\forall~  \dot\a\in \dot R\setminus\dot S(i)~  \exists\;\; N\in\bbbz^{>0}\ni ~~\fl^{\dot \a+n\d} v=\{0\} \hbox{(resp. $\fl^{\dot \a-n\d} v=\{0\}$)}\quad (\forall n\geq N).$
\end{itemize}
\end{Pro}
\pf As the proof of the case that $S(i)$ is down-nilpotent hybrid is similar to the proof of the case that $S(i)$ is up-nilpotent hybrid, we just give the  proof of the case that $S(i)$ is up-nilpotent hybrid. This is  a   modified version of  what we  give to prove (5.10) of \cite{you8}; but as the proof is so technical,  for the convenience of readers, we give the proof.

\noindent{$\bullet$} $\pmb{\fl\not=A(2k-1,2\ell-1)^{(2)}}:$
 By Proposition~\ref{submod}(iii), there is a functional $\boldsymbol{f}$ on $\sspan_\bbbr R$ such that $\boldsymbol{f}(\d)>0$ and
\begin{align*}
F:=&\{\lam\in\supp(M)\mid  (\lam+(S(i)\cap R^+))\cap \supp(M)=\emptyset\}\\
\cap & \{\lam\in\supp(M)\mid \forall\a\in \dot S(j) ~\exists N\in\bbbz^{>0}~\ni \lam+\dot\a+\bbbz^{\geq N}\d\cap \supp(M)=\emptyset\}
\end{align*}
is a  nonempty subset of $\supp(M)$.
 Set
\begin{align*}
\aa:=&\{{0\neq v}\in M\mid    \fl^{\a} v=\{0\}\quad(\a\in S(i)\cap R^+)\}\\
\cap & \{v\in M\mid \forall  \dot\a\in {\dot R_{re}\setminus \dot S(i)=\dot S(j)_{re}^\times}\;\;\exists N\in \bbbz^{>0}\ni \fl^{\dot \a+n\d} v=\{0\}\;\;\;(\forall n\geq N)\}.
\end{align*}
Pick $\lam\in F$  and fix $0\neq v\in M^\lam.$ Then, $v\in\aa.$ So to complete the proof, we just need to show that for each $\dot\a\in \dot R_{ns}^\times,$ there is a positive integer $N$ such that for each $n\geq N,$ $\fl^{\dot\a+n\d}v=\{0\}.$

Suppose that $\dot\a\in \dot R_{ns}^\times.$ Then by Table~\ref{table1}, we have
 $\dot\a=\dot\b+\dot\gamma,$ for some $\dot\b,\dot\gamma\in \dot R_{sh},$ and by (\ref{imp}(b)), we have
$$ S_{\dot\b}=S_{\dot\gamma}=\bbbz\d\hbox{~~as well as ~~} S_{\dot\a}=s\bbbz\d\quad\quad(\hbox{for some $s\in\bbbz^{>0}$}).$$
This together with the fact that  $\boldsymbol{f}(\d)>0$ and  $\lam\in F,$ guarantees the existence  of a  large enough $n$ such that $\lam+\dot\b+sn'\d,\lam+\dot\gamma+sn'\d\not\in \supp(M)$ for all $n'\geq n.$
So, for each  nonnegative integer $t$, we have
\begin{align*}
\fl^{ \dot\a+s(2n+t)\d}v=&[\fl^{ \dot\b+s(n+t)\d},\fl^{ \dot\gamma+sn\d}]v\sub \fl^{ \dot\b+s(n+t)\d}{M^{\lam+\dot\gamma+sn\d}}+ \fl^{ \dot\gamma+sn\d}{M^{\lam+\dot\b+s(n+t)\d}}=\{0\}.
\end{align*}
Therefore, for each $n'\geq 2n,$ we have $\fl^{ \dot\a+sn'\d}v=\{0\}$ and so we are done as $S_{\dot\a}=s\bbbz\d$.

\noindent{$\bullet$} $\pmb{\fl=A(2k-1,2\ell-1)^{(2)}}:$ We have $R_{re}\sub R_0.$ Recalling \cite[Rem.~3.1]{you8}, if $\a\in R_{re}\cup R_{im}$  and  $\ep\in R_{ns}^\times$  with $\a+\ep\in R,$ then, we have $\ep+\a\in R_{ns}$.
Set $$W:=\sum_{\lam\in \supp(M)}\sum_{\ep\in R_{ns}^\times}\fl^\ep M^\lam.$$
For $\a\in S(i)$ $(i=1,2),$ we have
\begin{align*}
\fl^\a W&=\fl^\a\hbox{\small$\displaystyle{\sum_{\lam\in \hbox{\tiny{\rm supp}}(M)}\sum_{\ep\in R_{ns}^\times}}$}\fl^\ep M^\lam=\hbox{\small$\displaystyle{\sum_{\lam\in \hbox{\tiny{\rm supp}}(M)}\sum_{\ep\in R_{ns}^\times}}$}\fl^\a\fl^\ep M^\lam\\
&\sub \hbox{\small$\displaystyle{\sum_{\lam\in \hbox{\tiny{\rm supp}}(M)}\sum_{\ep\in R_{ns}^\times}}$}\underbrace{[\fl^\a,\fl^\ep]}_{\in \sum_{\eta\in R_{ns}^\times}\fl^\eta} M^\lam+ \hbox{\small$\displaystyle{\sum_{\lam\in \hbox{\tiny{\rm supp}}(M)}\sum_{\ep\in R_{ns}^\times}}$}\fl^\ep\underbrace{\fl^\a M^\lam}_{\in\sum_{\mu\in \hbox{\tiny{\rm supp}}(M)}M^\mu}\sub W;
\end{align*}
in other words, $W$ is a $\scg(i)$-submodule of $M.$
Using   Proposition~\ref{submod}(iii),
\begin{equation}\label{new9}
\parbox{4.3in}{\small there are a fuctional $\bs{f}$ on $\sspan_\bbbr R$ with corresponding triangular decomposition $R=R^+\cup R^\circ\cup R^-$ such that $\bs{f}(\d)>0$ and  a weight $\mu $ of $W$ such that
$\mu+\a$ is not a weight of $W$ if $\a\in  S(i)\cap R^+$ and moreover, for each $\dot\a\in \dot S(j),$ there is a positive integer  $N$ such that $\mu+\dot \a+n\d\not\in \supp(W)$ for all $n\geq N.$}
\end{equation}
 Since $\mu$ is a weight of $W,$  there is a nonzero nonsingular root ${\ep_*}$ and $\lam\in \supp(M)$ such that $\fl^{\ep_*} M^\lam\neq \{0\}$ and $\mu={\ep_*}+\lam.$ For $0\neq v\in \fl^{\ep_*} M^\lam,$ we have
\begin{equation}\label{not-root-w}\fl^\a v\sub W^{\a+\mu}\stackrel{(\ref{new9})}{=\joinrel=}\{0\}\quad\quad (\a\in R^+\cap S(i));
\end{equation}
and
\begin{equation}\label{not-root-w}
\hbox{ $\forall\dot\a\in \dot S(j)$  $\exists N\in \bbbz^{\geq 0}$ such that $\fl^{\dot \a+n\d} v\sub W^{\dot\a+n\d+\mu}\stackrel{(\ref{new9})}{=\joinrel=}\{0\}\quad (n\geq N)$}.
\end{equation}
To complete the proof, we need to show (\ref{not-root-w}) holds for  all $\dot\a\in \dot R\setminus \dot S(i)^\times=\dot R_{ns}^\times \cup \dot S(j).$ We first note that $\dim(\fl^{\ep_*})=1$ and that two times of a {nonzero} nonsingular root is not a root (in particular, $\fl^{2{\ep_*}}=\{0\},$) so
\begin{equation}
\label{g-ep}
\fl^{\ep_*} v\sub\fl^{\ep_*}\fl^{\ep_*} M^\lam\sub{[\fl^{\ep_*},\fl^{\ep_*}]}M^\lam=\{0\}.
\end{equation}
Suppose  $${\ep_*}={\dot\ep_*}+s\d\quad \hbox{for some ${\dot\ep_*}\in \dot R^\times_{ns} $ and $s\in \bbbz.$}$$ For each $\dot\a\in \dot R_{ns}^\times,$
by \cite[Rem.~3.1]{you8}, one of the following happens:
\begin{itemize}
\item $\exists\;\dot\b_1\in \dot R_{sh}\ni \dot \a={\dot\ep_*}+\dot \b_1,$
    \item $\exists \; \dot\b_1\in \dot R_{sh}, \dot\b_2\in \dot R_{re}^\times \ni
     {\dot\ep_*}+\dot \b_1\in \dot R_{ns}^\times,
     \dot \a={\dot\ep_*}+\dot \b_1+\dot\b_2,$
     \item $\exists\;\dot\b_1\in \dot R_{sh},\dot \b_2,\dot\b_3\in \dot R_{re}^\times\ni
     {\dot\ep_*}+\dot \b_1,{\dot\ep_*}+\dot\b_1+\dot\b_2\in \dot R_{ns}^\times,
     \dot \a={\dot\ep_*}+\dot \b_1+\dot\b_2+\dot\b_3.$
\end{itemize}

$\bullet$ In the first  case, using (\ref{new9}) together with (\ref{imp}(b)) and the fact that $\boldsymbol{f}(\d)>0$, we choose  $t_1\in \bbbz^{>0}$ such that
\begin{equation}\label{new00}
(\lam+\dot\b_1+\bbbz^{\geq t_1}\d)\cap  \supp(M)=\emptyset.
\end{equation} So, for $t>t_1+s,$  we have
\begin{align*}
\fl^{\dot \a+t\d}v=[\fl^{\dot\b_1+(t-s)\d},\fl^{ {\ep_*}}]v&\sub \fl^{\dot\b_1+(t-s)\d}\fl^{ {\ep_*}}v+\fl^{ {\ep_*}}\fl^{\dot\b_1+(t-s)\d}v\stackrel{(\ref{new00}),(\ref{g-ep})}{=\joinrel=\joinrel=\joinrel=}\{0\}.
\end{align*}

$\bullet$ In the second case, contemplating (\ref{imp}(b)), we use  (\ref{new9}) and the fact that $\boldsymbol{f}(\d)>0$  to choose
 $t_1,t_2\in\bbbz^{>0}$ with
\begin{equation*}\label{new01}
(\lam+\dot\b_1+\bbbz^{\geq t_1}\d)\cap\supp(M)=\emptyset,\;\; \dot\b_2+ t_2\d\in R\andd \lam+\dot\b_2+ t_2\d\not \in\supp(M).
\end{equation*}
 This implies that  for $t\geq  t_1+t_2+s,$ we have  $\fl^{\dot\b_1+(t-t_2-s)\d}v=\{0\}$ and $\fl^{\dot\b_2+t_2\d}v=\{0\}.$  So (\ref{g-ep}) implies that
$$\fl^{\dot \a+t\d}v=[\fl^{\dot\b_2+t_2\d},[\fl^{\dot\b_1+(t-t_2-s)\d},\fl^{ {\ep_*}}]]v=\{0\}.$$

$\bullet$
In the third case,  we
choose $t_1,t_2, t_3\in\bbbz^{>0}$ with
\begin{align*}
\hbox{\small $
(\lam+\dot\b_1+\bbbz^{\geq t_1}\d)\cap\supp(M)=\emptyset, ~\dot\b_i+ t_i\d\in R\andd \lam+\dot\b_i+ t_i\d\not \in\supp(M)\;\;\;(i=2,3).$}
\end{align*} Then for $t\geq  t_1+t_2+t_3+s,$  we have
$$\fl^{\dot \a+t\d}v=[\fl^{\dot\b_3+t_3\d},[\fl^{\dot\b_2+t_2\d},[\fl^{\dot\b_1+(t-t_2-t_3-s)\d},\fl^{ {\ep_*}}]]]v=\{0\}.$$
This completes the proof.
 \qed

\begin{Thm}\label{main}
Suppose that $\{r,t\}=\{1,2\}$ and  $M$ is a $t$-quasi-integrable $\fl$-module of nonzero level. Then there are triangular decompositions $R=R^+\cup R^\circ\cup R^-$ and $R^\circ=R^{\circ,+}\cup R^{\circ,\circ}\cup R^{\circ,-},$ with $\d\not\in R^\circ$  as well as $R^{\circ,\circ}\sub S(r),$ and a cuspidal module $\Omega$ over $\fl^\circ_{_P}$, for $P:=R^+\cup R^{\circ,+}\cup R^{\circ,\circ}$, such that $M\simeq{\rm Ind}_{_P}(\Omega).$ We also have that $\fl^\circ_{_P}$ is a direct sum of a reductive finite dimensional  Lie algebra and finitely many  basic classical simple Lie superalgebras  of types $B(0,p)$ ($p\in\bbbz^{>0}$). In particular, the classification of quasi-integrable irreducible finite weight $\fl$-modules is reduced to the classification of cuspidal modules over $\fl^\circ_{_P}$; see  \cite{M}, \cite{FGG} and \cite{G}.
\end{Thm}
\pf We assume $S(r)$ is up-nilpotent hybrid and  carry out the proof; the proof in case  $S(r)$ is down-nilpotent hybrid, is similarly done.

By Proposition~\ref{gen}, there is  a functional $\boldsymbol{f}$ on $\sspan_\bbbr R$ with $\boldsymbol{f}(\dot S(t))=\{0\}$ and $\boldsymbol{f}(\d)>0$ such that the set $\aa$ consisting of all nonzero weight vectors $v$ satisfying
\begin{itemize}
\item $\fl^{\a}v=\{0\}$ for all $\a\in S(r)$ with $\boldsymbol{f}(\a)>0,$
\item  $\forall\dot\a\in \dot R\setminus \dot S(r)~\exists ~m_{\dot\a}\in \bbbz^{>0}$ with $\fl^{\dot\a+k\d}v=\{0\}$ for all $k>m_{\dot\a}$
\end{itemize}
is nonempty.
For $v\in\aa,$ set
\[C_v:=\{\a\in R\setminus S(r)\mid \fl^{\a}v\neq \{0\},~\boldsymbol{f}(\a)>0\}.\] We mention that as $\boldsymbol{f}(\d)>0,$ $C_v$ is a finite set and claim that  there is $w_0\in \aa$ with $C_{w_0}=\emptyset.$ Pick $v\in \aa$ such that $C_v$ is of minimal cardinality. If $C_v=\emptyset,$ we take $w_0:=v$ and we are done. Otherwise, we pick $\a_*(v)\in C_v$ with $$\boldsymbol{f}(\a_*(v))={\rm max}\{\boldsymbol{f}(\a)\mid \a\in C_v\}.$$ Suppose $0\neq w_1\in \fl^{\a_*(v)}v.$ We claim that $$w_1\in \aa\andd C_{w_1}= C_v.$$

\noindent $\bullet ~{w_1\in \aa}:$ We show it in the following two steps:

\noindent \underline{Step~1}. For $\a\in S(r)$ with $\boldsymbol{f}(\a)>0,$ since $v\in\aa,$ we have
\begin{align*}
\fl^{\a}w_1\sub\fl^\a\fl^{\a_*(v)}v\sub [\fl^{\a},\fl^{\a_*(v)}]v+\fl^{\a_*(v)}\underbrace{\fl^\a v}_{\{0\}}\sub \fl^{\a+\a_*(v)}v.
\end{align*}
We have   \[\boldsymbol{f}(\a+\a_*(v))>\boldsymbol{f}(\a_*(v))>0.\] If $\a+\a_*(v)\not\in R$ or $\a+\a_*(v)\in S(r),$ we get $\fl^{\a+\a_*(v)}v=\{0\}$ as $v\in\aa$ and so $\fl^{\a}w_1\sub \fl^{\a+\a_*(v)}v=\{0\}.$ Also if $\a+\a_*(v)\in R\setminus S(r),$ due to the choice of $\a_*(v)$, $\a+\a_*(v)\not\in C_v$ and so  $\fl^{\a+\a_*(v)}v=\{0\}.$  Therefore, again we have  $\fl^{\a}w_1\sub \fl^{\a+\a_*(v)}v=\{0\}.$

\noindent \underline{Step~2}. Since $v\in \aa,$ we choose $N$ such that
\begin{equation}\label{cv}
\boldsymbol{f}(\dot\a+n\d)>0\andd \fl^{\dot\a+n\d}v=\{0\}\quad\quad (\dot\a\in \dot R,~n>N).
\end{equation}
Suppose $\dot\a\in \dot R\setminus \dot S(r)$ and $n>N.$ If $\dot\a+n\d+\a_*(v)\in S(r),$ since $v\in\aa$ and $\boldsymbol{f}(\dot\a+n\d+\a_*(v))>\boldsymbol{f}(\a_*(v))>0,$ we get $\fl^{\dot\a+\a_*(v)+n\d} v=\{0\}.$ Also, if  $\dot\a+n\d+\a_*(v)\in R\setminus S(r),$  then as $\boldsymbol{f}(\dot\a+n\d+\a_*(v))>\boldsymbol{f}(\a_*(v)),$ we get  $\dot\a+n\d+\a_*(v)\not\in C_v$ and so again $\fl^{\dot\a+\a_*(v)+n\d} v=\{0\}.$ These altogether imply that for all $n> N,$ we have
\begin{align*}
\fl^{\dot\a+n\d}w_1\sub\fl^{\dot\a+n\d}\fl^{\a_*(v)}v\sub [\fl^{\dot\a+n\d},\fl^{\a_*(v)}]v+\fl^{\a_*(v)}\underbrace{\fl^{\dot\a+n\d} v}_{(\ref{cv})}\sub \fl^{\dot\a+n\d+\a_*(v)}v=\{0\}
\end{align*}
as we desired.
\smallskip

\noindent$\bullet~{C_{w_1}= C_v:}$  Suppose that $\a\in C_{w_1}.$ Since $\boldsymbol{f}(\a+\a_*(v))>\boldsymbol{f}(\a_*(v))>0,$ we get $\fl^{\a+\a_*(v)}v=\{0\}.$
So,  we have
\begin{align*}
\{0\}\neq \fl^{\a}w_1\sub \fl^\a\fl^{\a_*(v)}v\sub [\fl^{\a},\fl^{\a_*(v)}]v+\fl^{\a_*(v)}\fl^\a v\sub \underbrace{\fl^{\a+\a_*(v)}v}_{\{0\}}+\fl^{\a_*(v)}\fl^\a v
\end{align*}
which in turn implies that  $\fl^\a v\neq\{0\},$ that is, $\a\in C_v.$ So $C_{w_1}\sub C_v.$
Since $C_v$ is of minimal cardinality, we get that $C_v=C_{w_1},$ as we desired.

Since $C_{w_1}=C_v$, we have in particular that   $$\a_*(v)\in C_{w_1}\andd \boldsymbol{f}(\a_*(v))={\rm max}\{\boldsymbol{f}(\a)\mid \a\in C_{w_1}\}.$$ As $\a_*(v)\in C_{w_1}$  and $\fl^{\a_*(v)}$ is 1-dimensional, we get $\{0\}\neq \fl^{\a_*(v)}w_1= \fl^{\a_*(v)}\fl^{\a_*(v)}v$ and $\a_*(v)\in R_{ns}\cup S(t).$ If $\a_*(v)\in R_{ns},$ since $R_{ns}\sub R_1$ and  two times of a nonzero nonsingular root is not a root, we get $\{0\}=[\fl^{\a_*(v)},\fl^{\a_*(v)}]v=\fl^{\a_*(v)}\fl^{\a_*(v)}v$ which is a contradiction. So
\begin{equation}\label{sj}
\a_*(v)\in S(t)_{re}^\times.
\end{equation}
Repeating  the above process for $w_1$ instead of $v,$ we get $0\neq w_2\in \fl^{\a_*(v)}w_1\sub  \fl^{\a_*(v)}\fl^{\a_*(v)}w_1$ with
  $$C_v=C_{w_1}=C_{w_2},~~\a_*(v)\in C_{w_1}=C_{w_2}~\hbox{ and }~ \boldsymbol{f}(\a_*(v))={\rm max}\{\boldsymbol{f}(\a)\mid \a\in C_{w_1}=C_{w_2}\}.$$
Continuing this process,    for each $n\in\bbbz^{>0},$
\[\{0\}\neq\underbrace{\fl^{\a_*(v)}\cdots\fl^{\a_*(v)}}_{n~{\rm times}}v\sub \fl^{\lam+n\a_*(v)}\] which is a contradiction as by our assumption and (\ref{sj}), $\a_*(v)\in R^{ln}\sub  \frak{B}_M$. So, $C_v$ cannot be nonempty. This means that $\fl^+v=\{0\}$. So by \cite[Pro.~3.3]{you8},
 $N:=\{w\in V\mid \fl^+w=\{0\}\}$ is an irreducible module over $\op_{\a\in R^\circ}\fl^\a$ and for $P':=R^+\cup R^\circ,$ $M\simeq {\rm Ind}_{_{P'}}(N).$ {Since $\d\not \in R^\circ,$ $R^\circ$ is finite and   $\op_{\a\in R^\circ}\fl^\a$ is a finite dimensional Lie superalgebra.
 Since \[{\rm inj}N:=R^\circ\cap R^{in}\cap R_0\sub R(r),\] using \cite[Thm.~3.6]{DMP}, one has a triangular decomposition $R^\circ=R^{\circ,+}\cup R^{\circ,\circ}\cup R^{\circ,-}$ for $R^\circ$ with  $R^{\circ,\circ}\sub S(r)$ and a cuspidal module $\Omega$ over $\sum_{\a\in R^{\circ,\circ}}\fl^\a$ such that for  $P'':=R^{\circ,+}\cup R^{\circ,\circ},$  $N\simeq{\rm Ind}_{_{P''}}(\Omega).$ Therefore, \[M\simeq {\rm Ind}_{_{P'}}(N)\simeq {\rm Ind}_{_{P'}}({\rm Ind}_{_{P''}}(\Omega))\simeq{\rm Ind}_{_P}(\Omega),\] where $P=R^{+}\cup R^{\circ,+}\cup R^{\circ,\circ}.$ As $R^{\circ,\circ}\sub S(r)$ just contains real roots, it is a direct sum of a finite root system and finitely many root systems of types $B(0,p)$ ($p\in\bbbz^{>0}$)} and so we are done.
\qed

\bibliographystyle{amsplain}

\end{document}